\documentclass[a4paper,12pt]{article}
\usepackage{graphicx}
\usepackage{psfrag}

\usepackage{mathrsfs}
\usepackage{amssymb}
\usepackage{amsfonts}
\usepackage{latexsym}
\usepackage{amsmath}
\usepackage{amsthm}
\usepackage[cp1251]{inputenc}
\usepackage[english]{babel}

\newcommand{\er}[1]{\textrm{(\ref{#1})}}
\def\lb{\label}

\theoremstyle{plain}

\newtheorem{theorem}{\bf Theorem}[section]
\newtheorem{lemma}[theorem]{\bf Lemma}

\newtheorem{corollary}[theorem]{\bf Corollary}
\theoremstyle{remark}

\renewcommand{\b}{\beta}                 
\newcommand{\g}{\gamma}                  
\newcommand{\G}{\Gamma}

\newcommand{\D}{\Delta}                
\newcommand{\ve}{\varepsilon}

\newcommand{\vt}{\vartheta}       \newcommand{\cJ}{\mathcal{J}}     
          \newcommand{\cK}{\mathcal{K}}     
                
\renewcommand{\l}{\lambda}               
               
\newcommand{\m}{\mu}                   
\newcommand{\n}{\nu}                   
\renewcommand{\r}{\rho}                  
\newcommand{\s}{\sigma}

\newcommand{\vp}{\varphi}              
                   
\newcommand{\p}{\psi}                   
             \newcommand{\cZ}{\mathcal{Z}}

\def\Z{\mathbb{Z}}
\def\R{\mathbb{R}}
\def\C{\mathbb{C}}
\def\T{\mathbb{T}}

\def\K{\mathbb{K}}

\def\qqq{\qquad}
\def\qq{\quad}
\let\ge\geqslant
\let\le\leqslant

\newcommand{\ca}{\begin{cases}}
\newcommand{\ac}{\end{cases}}
\newcommand{\ma}{\begin{pmatrix}}
\newcommand{\am}{\end{pmatrix}}

\def\lt{\biggl}
\def\rt{\biggr}

\def\eq{\begin{equation}}
\def\qe{\end{equation}}
\renewcommand{\[}{\begin{equation}}
\renewcommand{\]}{\end{equation}}

\def\pa{\partial}
\def\sm{\setminus}
\def\es{\emptyset}
\def\no{\noindent}
\def\ol{\overline}
\def\iy{\infty}
\def\ev{\equiv}
\def\/{\over}
\def\ts{\times}

\def\os{\oplus}

\def\Re{\mathop{\rm Re}\nolimits}

\def\Im{\mathop{\rm Im}\nolimits}

\def\BBox{\hspace{1mm}\vrule height6pt width5.5pt depth0pt \hspace{6pt}}

\begin{document}
\title{ A priori estimates for the Hill and Dirac operators}
\author{
Evgeny Korotyaev
\begin{footnote} {
Institut f\"ur  Mathematik,  Humboldt Universit\"at zu Berlin,
Rudower Chaussee 25, 12489, Berlin,Germany, e-mail:
evgeny@math.hu-berlin.de 
}
\end{footnote}
}
\maketitle

\begin{abstract}
\no Consider the Hill operator $Ty=-y''+q'(t)y$ in $L^2(\R)$, where 
$q\in L^2(0,1)$ is a 1-periodic real potential. 
The spectrum of $T$ is is absolutely
continuous and consists of bands  separated by gaps $\g_n,n\ge 1$ with length $|\g_n|\ge 0$. We obtain a priori estimates of the gap lengths,
effective masses, action variables for the KDV. 
For example, if $\m_n^\pm$ are the effective masses associated with
the gap $\g_n=(\l_n^-,\l_n^+)$, then $|\m_n^-+\m_n^+|\le C|\g_n|^2n^{-4}$
for some constant $C=C(q)$ and any $n\ge 1$.
In order prove these results we use the analysis of a
conformal mapping corresponding to quasimomentum of the Hill operator. That makes possible to reformulate the problems for the differential operator as the problems of the conformal mapping theory. Then the proof is based on the analysis of the conformal mapping  and the identities.
Moreover, we  obtain the similar  estimates for the Dirac operator.
\end{abstract}

\section{Introduction  and main results}
\setcounter{equation}{0}

Consider the self-adjoint periodic weighted operator  $T_w$  in $L^2(\R,w^2(t)dt)$, given by
$$
T_wf=-w^{-2}(w^2f')'=-f''-2pf',  \ \ w(t)=e^{\int_0^xp(s)ds}, \ \ p\in H=\{p\in L^2(\T), \int_0^1 p(t)dt=0\}.
$$
The operator $T_w$ is well studied, see \cite{K8}, \cite{K9}, \cite{Kr}
 and references therein.
Define the unitary transformation $U:L^2(\R,w^2 dt)\to L^2(\R,dt)$
as multiplication by $w$. Then $T_w$ is unitarily equivalent to
the Hill operator $T=UT_wU^{-1}$ given by
\[
\lb{1.10}
T=-\lt({d \/dt}+p\rt)\lt({d \/ dt}-p\rt)=
-{d^2\/dx^t}+q'+q_0\ge 0,\ \ \ \ \ q'=p'(t)+p^2(t)-\|p\|^2,
\]
where  $q_0=\|p\|^2=\int_0^1 p^2(t)dt$ and 
$q'$ is  a 1-periodic potential (distribution). We mention following papers devoted to singular potentials 
\cite{KM}, \cite{K7}, \cite{NS}, \cite{SS}. We define the Riccati type mapping $R:H\to H$ by 
\[
\lb{1.11}
q=R(p)=p(x)+\int_0^x(p^2(t)-\|p\|^2)dt+\int_0^1\rt(t-{1\/2}\rt)p^2(t)dt,\ \ \ p\in H.
\]
Recall that the Riccati map $R:H\to H$ is a real analytic isomorphism
between $H$ and $H$, see \cite{K7}. Thus, if $p'\in L^2(\T),$ then $T_w$ corresponds to the Hill  operator $T$ with $L^2$-potential.

We can not introduce the standard fundamental solutions for the operator $T$, since the perturbation $q'$ is very strong. But 
the operator $T_w=U^{-1} T U$ has the standard fundamental solutions $\vp(t,\l), \vt(t,\l)$, which satisfy the  equation
$-y''-2py'=\l y, \ \l\in \C$ and the conditions
$\vp'(0,\l)=\vt(0,\l)=1,\ \vp(0,\l)=\vt'(0,\l)=0$.
Here and below we use the notation $u'={\pa u\/\pa t}$.
Introduce the Lyapunov function $\D(\l)={\vp'(1,\l)+\vt(1,\l)\/2}$.

It is well known, that the spectrum $\s(T)=\s(T_w)$ of $T$  is absolutely continuous and consists of intervals $\s _n=[\l^+_{n-1},\l^-_n ],$  where
$\l^+_{n-1}<\l^-_n \le  \l^+_{n},\ n\ge 1$ and  
$\l_0^+=0$.
The intervals $\s_{n}$ and
$\s_{n+1}$  are separated by the gap $\g_n=(\l^-_n,\l^+_n )$.
 If a gap degenerates, that is $\g_n=\es$,
then the corresponding segments $\s_{n} $ and $\s_{n+1}$ merge.
 Note that $\D(\l_{n}^{\pm})=(-1)^n,\ n\ge 1,$ and
there exists a unique point  $\l_n\in [\l^-_n,\l^+_n]$ such that
 $\D'(\l_n)=0$. The sequence $\l_0^+<\l_1^-\le \l_1^+\ <\dots$
 is the spectrum of equation $-y''+qy=z y$ with the 2-periodic
 boundary conditions, i.e. $y(t+2)=y(t), t\in \R$.
If $\l_n^-=\l_n^+$ for some $n$, then this number $\l_n^{\pm}$
 is the double eigenvalue of this equation with the 2-periodic
 boundary conditions.
The lowest  eigenvalue $\l_0^+$ is always simple, $\D(\l_0^+)=1,$  and the
corresponding eigenfunction is 1-periodic. The eigenfunctions,
corresponding to the eigenvalue $\l_n^{\pm}$, are 1-periodic,
when $n$ is even and they are antiperiodic,  i.e. $y(t+1)=-y(t),\ \
t\in\R$, when $n$ is odd. 

We introduce a new parameter (momentum) $z=x+iy=\sqrt\l, \sqrt 1=1$,  and  momentum gaps
$$
g_n=(z^-_n,z^+_n ),\qq z_n^{\pm}=\sqrt{\l_n^{\pm}}\ge 0,\qqq g_{-n}=-g_n,\qq n\ge 1,\qq g_0=\es.
$$
We define a quasimomentum  
$k(z)=u(z)+iv(z)=\arccos \D(z^2), \ \ z\in \cZ={\bf C}\sm \ol{g}, \ g=\cup g_n$, see [MO].
The function $k(z)$ is  a conformal mapping from $\cZ$ onto a quasimomentum domain $\cK=\C\sm\cup_{n\in \Z}\G_n$, where the slit 
$\G_n=[\pi n+ih_n,\pi n-ih_n],\ h_0=0$,
and $h_n=h_{-n}\ge 0$ is defined by the equation $\D(\l_n)=(-1)^n\cosh h_n$. The quasimomentum $k$ satisfies $k(z)=z+o(1)$ as $|z|\to \iy$.
With  each edge  of  the gap $\g_n\ne \es$, we associate the effective mass $\m_0^+, \m_n^{\pm}$ by  
\[
\l=z^2=\l_n^{\pm}+{(k-\pi n)^2\/2\m_n^{\pm}}(1 +o(1))
\ \ \ \ \  {\rm as} \ \ \ \             \l\to \l_n^{\pm}.
\]
Introduce the real Hilbert spaces 
$\ell^2_p=\rt\{f=(f_n)_1^\iy ,\ \ f_n\in \ \R,\ 
\| f\|_p^2=\sum_{n\ge 1}(2\pi n)^{2p}f_n^2<\iy \rt\},\ \ p\ge 0$ and let $\ell_0^2\ev\ell^2$. 
For the sake of the reader, we briefly recall the results existing in the literature about the a priori estimates.
Marchenko and Ostrovki [MO1-2] obtained  the estimates: 
$\|q_1\|\le C(1+h_{+})\|h\|_1,\ \ \|h\|_1\le C\|q_1\|\exp (\ A\|q_1\|\ ),$
for some absolute constants $A, C$, where $h_+=\sup_{n\ge 1}h_n$
and $q_1=q_0+q'$. These estimates are not sharp
since they used the Bernstein inequality.
Using the harmonic measure argument  Garnett and Trubowitz \cite{GT} obtained $\|\g\|\le (4+\|h\|_1)\|h\|_1$. 

Various estimates for the Hill
operator with a potential $q'\in L^2(0,1)$ were obtained in the author's paper [K5-6] and in [KK1], [K3] for the Dirac operator. 
In the paper [K3], [K5], [K6] the following sharp estimates for the case $q'\in L^2(0,1)$ were obtained:
$$
\|\g \|\le 6\|q'\|( 1+\|q'\|)^{1\/3} ,\ \ \ \ \ 
\|q'\|\le 4\|\g\|( 1+\|\g \|)^{1\/3} , \qqq \ \ \g=(|\g_n|)_1^\iy,
$$$$
2\|h\|_1\le \pi \|q'\|( 1+\|q'\|)^{1\/3} ,\ \ \ \ 
\|q'\|\le 3\|h\|_1(6+h_{+})^{1\/2},\ \ \ \ h_+=\sup_{n\ge 1}h_n,
$$
and 
 for the case $q\in L^2(0,1)$ [K8]:
$$
\|q\|\le 2\|h\|(1+4\|h\|),\ \ \ \ \|h\|\le 3\|q\|(1+2\|q\|),
\ \ \ \ q\in H,
$$
$$
\|q\|\le 48\pi^2 \|\g\|_{-1}(1+\|\g\|_{-1})^3,\ \ \ \
\|\g\|_{-1}\le 2 \|q\|(1+2\|q\|)^3.
$$

It is well known that $v(z)=\Im k(z)\ge v_n(z)=|(z-z_n^-)(z-z_n^+)|^{1\/2}, z\in g_n$, see [KK1]. Define a function $V_n(z)={v(z)\/v_n(z)}-1\ge 0, z\in g_n$, which is analytic in the domain $\C_+\cup \C_-\cup (z_{n-1}^+,z_{n+1}^-)$, see [KK1]. Let 
$M_n=\max_{z\in g_n} V_n(z), M_n'=\max_{z\in g_n} V_n'(z),..$

\begin{theorem} \label{T1} 
Let real $q\in L^2(0,1)$. 
Then for each $n\ge 1$ the following estimates hold
\[
\lb{T1-1} 
|\l_n^0-{|g_n|^2\/4}-\l_n| \le {3|\g_n|^2\/8z_n^0}M_n',\qq where
\qq z_n^0={z_n^-+z_n^+\/2},\qq \l_n^0={\l_n^-+\l_n^+\/2},
\]
\[
\lb{T1-2}
|h_n-{|\g_n|(1+V_n(z_n))\/4z_n^0}|\le {|\g_n|^3M_n'^2\/16{z_n^0}^3},
\]
\[
\lb{T1-3}
|h_n-h_n^{\pm }|\le{|g_n|^2\/2}M'_n\rt(1+
{M_n'\sqrt{|\g_n||\m_n^\pm|}\/4}\rt), \qq where\ 
h_n^{\pm }=\rt||g_n|z_n^{\pm}\m_n^{\pm}\rt|^{1\/2},
\]
\[
\lb{T1-4}
|\m_n^++\m_n^-|\le {|g_n|^2\/4z_n^+z_n^-}(1+M_n)(1+2z_n^0M_n'),
\]
\[
\lb{T1-5}
|\m_n^{\pm}\mp \m_n|\le 
{|g_n|^2\/8z_n^0z_n^\pm}
(1+M_n)^2\rt[1+4z_n^\pm M_n'\rt(1+{M'_n\/4}\sqrt{|\g_n||\m_n^\pm|}\rt)\rt],\ \ where \ \m_n={2h_n^2\/|\g_n|}.
\]
\end{theorem} 

\no {\bf Remark.} 1) Using Theorem \ref{T1} and estimates
\[
|\g_n|\le (4\pi n)^2(\m_n^+-\m_n^-),\qqq
h_n^2\le {(3\pi)^2\/2}|\g_n||\m_n^\pm|,
\]
from \cite{K4}, we obtain the estimates in terms of effective
masses.

2) Recall that $z_n^\pm=\pi n\pm
|\int _0^1e^{i2\pi nx}q(x)dx|+\ell^d(n),\ d>1$ as $n\to \iy$, see \cite{K8}. This yields asymptotic estimate $\r_n\ge \pi -\ve$
as $n\to \iy$, for any $\ve>0$.

Define the action variables $A_n
={4\/\pi }\int_{g_n}xv(x)dx\ge 0,\ n\ge 1$ for the KDV equation \cite{FM}. 

\begin{theorem} \label{T2} 
Let real $q\in L^2(0,1)$ and 
let $\b_n=\sqrt{A_n}\ge 0, n\ge 1$. Then the following estimates hold
\[
\lb{T2-1}
\rt|A_n-{|\g_n|^2\/8z_n^0}(1+V_n(z_n^0))\rt|\le {|g_n|^4\/2^5}(z_n^0 M_n''+M_n'),
\]
\[
\lb{T2-2}
|A_n-{|\g_n|h_n\/2}|\le {|g_n|^4\/2^5}(z_n^0M_n''+6z_n^0{M_n'}^2+M_n'),
\]
\[
\lb{T2-3}
|A_n-{|\g_n|\/2}\rt|{|g_n|m_n^\pm\/2}\rt|^{1\/2}|
\le {|g_n|^3\/4} \rt(z_n^0M_n'+{|g_n|\/8}M_n'+|g_n|M_n''\rt),
\]
\[
\lb{T2-4}
\rt|\b_n-|\g_n|\sqrt{(1+V_n(z_n^0))\/8z_n^0}\rt|\le {|g_n|^3\/2^5\sqrt{2z_n^0}}(z_n^0 M_n''+M_n'),
\]
\[
\lb{T2-5}
\rt|\b_n-\sqrt{|\g_n|h_n\/2}\rt|\le
{|g_n|^3\/2^5\sqrt{2z_n^0}}(z_n^0M_n''+6z_n^0{M_n'}^2+M_n'),
\]
\[
\lb{T2-6}
|\b_n-\sqrt{|\g_n|\/2}{||g_n|m_n^\pm|^{1\/4}\/2^{1\/4}}|\le 
{|g_n|^2\/4\sqrt{2z_n^0}} \rt(z_n^0M_n'+{|g_n|\/8}M_n'+|g_n|M_n''\rt).
\]
\end{theorem} 

\no {\bf Consider the Dirac operator} $T_D$ acting in
the space $L^2(\R )\os L^2(\R )$ and given by
$$
T_D=\cJ {d\/ dt}+q_D(t),\ \ \ \ q_D=\ma q_1&q_2\\q_2&-q_1\am
\ \ \ \cJ\ev \ma 0&1\\-1&0\am ,
 $$
where  $q_1,q_2\in L^2(0,1)$ are real 1-periodic functions  of  $t\in \R$. The spectrum of $T_D$ is purely absolutely continuous and  $\s(T_D)=\cup_{n\in \Z} s_n$, where
$s_n=[z_{n-1}^+,z_n^-]$ are spectral bands and $z_{n-1}^+<z_n^-\le z_n^+$, (see [LS] for details). The intervals $s_n,s_{n+1}$ are separated by the gap $g_n=(z^-_n,z^+_n)$ with length $|g_n|\ge0$.
 If a gap $g_n=\es$,  then the corresponding
spectral bands $s_n, s_{n+1}$ merge. 
Introduce the $2\ts 2$-matrix valued fundamental solution
 $\p =\p (t,z)$  satisfying
$
J\p '+q_D\p=z \p,   \p (0,z)=I_2, \ \ z\in \C,
$
where $I_2$ is the $2\ts 2$ identity matrix.
Define the Lyapunov function $\D(z)={1 \/ 2}{\rm Tr} \p (1,z)$. 
Note that  $\D(z_n^{\pm})=(-1)^n, n\in\Z$, and the function
$\D'(z)$ has a zero $z_n$ in a ''closed gap'' $[z^-_n,z^+_n]$
(see e.g. [LS]). The sequence $..z_0^+<z_1^-\le z_1^+ <\dots$
 is the spectrum of equation $Jy'+q_Dy=z y, y=(y_1,y_2)^\top,$ with the 2-periodic  boundary conditions, i.e. $y(t+2)=y(t), t\in \R$.
If $z_n^-=z_n^+$ for some $n$, then $z_n^{\pm}$
 is the double eigenvalue of the equation $Jy'+q_{zs}y=z y$ with the 2-periodic  boundary conditions. The eigenfunctions,
corresponding to the eigenvalue $z_n^{\pm}$, are 1-periodic,
when $n$ is even and they are antiperiodic,  i.e. $y(t+1)=-y(t),\ \
t\in\R$, when $n$ is odd.

For each $q_D$ there exists a unique conformal mapping (the quasimomentum) $k:\cZ\to \cK=\C\sm\cup \G_n, \G_n=[\pi n-ih_n,\pi n+ih_n],$ given by (see \cite{K6})
$$
\cos k(z)= \D(z),\ \  z\in \cZ =\C\sm \ol g,\qq g=\cup  g_n,
\qq k(z)=z+o(1) \ as \ \ |z|\to \iy,
$$
where  the height $h_n\ge 0$ is defined by the equation $\cosh h_n = (-1)^n\D(z_n) \ge 1$. 
With  each edge  of  the gap $g_n\ne \es$, we associate the effective mass
$m_n^{\pm}$ by
\[
\lb{demzs}
z=z_n^{\pm}+{(k(z)-\pi n)^2 \/  2m_n^{\pm}}(1+o(1))
\qqq {\rm as} \ \ z\to z_n^{\pm}.
\]

\begin{theorem} \label{TD1} 
Let real $q_D\in L^2(0,1)$. Then for each $n\in \Z$ the following estimates hold
\[
\lb{TD1-1}
|z_n-z_n^0|\le {|g_n|^2\/4}M_n',
\]
\[
\lb{TD1-2}
|2h_n-|g_n|(1+V_n(z_n)|\le {|g_n|^3\/8} {M_n'}^2,
\]
\[
\lb{TD1-3}
|h_n-h_n^{\pm }|\le{|g_n|^2\/2}M_n'\rt(1+
{M_n'\/8}\sqrt{2|g_n||m_n^\pm|}\rt), \qq where\ 
h_n^{\pm }=\sqrt{|g_n||m_n^{\pm}|\/2}, 
\]
\[
\lb{TD1-4}
|m_n^++m_n^-|\le |g_n|^2(1+M_n)M_n',
\]
\[
\lb{TD1-5}
|m_n^{\pm}\mp m_n|\le |g_n|^2(1+M_n)M_n'
\rt(1+{M'_n\/8}\sqrt{2|g_n||m_n^\pm|}\rt), \qq where\ m_n={2h_n^2\/|g_n|}.
\]

\end{theorem}

{\bf Remark.} Using Theorem \ref{TD1} and estimates
\[
{|g_n|\/2}\le \ h_n\le \pi 
\sqrt {2|g_n||m_n^{\pm}|} \le 2\pi |m_n^{\pm}|,
\qq h_n^2\le 2|g_n| \sqrt {m_n^+|m_n^-|},\ \
 \ \ |g_n|\le 2|m_n^{\pm}|.
\]
from \cite{KK1}, we obtain the estimates in terms of effective
masses.

Recall that  the action variables for the NLS equation are given by
$a_n={1\/\pi }\int_{g_n}v(x)dx, \ n\in \Z$.

\begin{theorem} \label{TD2} 
Let real $q_D\in L^2(0,1)$ and let $b_n=\sqrt{a_n}\ge 0, n\in \Z$.
Then the following estimates hold
\[
\lb{TD2-1}
|a_n-(|g_n|^2/8)(1+V_n(z_n^0))|\le 2{|g_n|^4\/4^4} M_n'',
\]
\[
\lb{TD2-2}
|a_n-{|g_n|h_n\/4}|\le {|g_n|^4\/2^7}(M_n''+6{M_n'}^2),
\]
\[
\lb{TD2-5}
|a_n-{|g_n|^{3/2}\/2^{5\/2}}|m_n^\pm|^{1\/2}|\le {|g_n|^3\/16} M_n'+{|g_n|^4\/2^7} M_n'',
\]
\[
\lb{TD2-3}
|b_n-2^{-3/2}|g_n|\sqrt{1+V_n(z_n^0)}|\le {\sqrt 2|g_n|^{3}\/4^3} M_n'',
\]
\[
\lb{TD2-4}
|b_n-\sqrt{|g_n|h_n/2}|\le {\sqrt 2 |g_n|^3\/2^7}(M_n''+6{M_n'}^2),
\]
\[
\lb{TD2-6}
|b_n-{|g_n|^{3\/4}\/2^{5\/4}}|m_n^\pm|^{1\/4}|\le {({|g_n|^2\/8} M_n'+{|g_n|^3\/8^2} M_n'')\/\sqrt2}.
\]

\end{theorem} 

Note that Theorem \ref{TD2} is used to study electrostatic inverse problems  on the plane \cite{KK3}.

Define integrals $Q_p={1\/\pi}\int_\R x^{p}v(x)dx, p\in \Z$.
Recall the identities for the Hill operator 
\[
Q_0={q_0\/2}={\|p\|^2\/2}, 
\qqq Q_2={1\/8}\int_0^1(q_0^2+q'(t)^2)dt,
\]
(recall that $\l_0^+=0$ for the Hill operator), and for the Dirac operator 
\[
Q_0={\|q_D\|^2\/2} ={1\/2}\int_0^1(q_1^2(t)+q_2^2(t))dt,
\qqq Q_2={1\/
8}\int_0^1(({q_1'}^2+{q_2'}^2)+(q_1^2+q_2^2)^2)dt.
\]

\begin{theorem} \label{T5} 
Let real $q\in L^2(0,1)$ or $q_D\in L^2(0,1)$. 
 Let $\r_n={\rm dist}(g_n, g\sm g_n),\  n\in \Z$.
Then for each $n\in \Z$ the following estimates hold
\[
\lb{T5-1}
 M_n \le {2\/\pi \r_n}\sup_{s\in \Z} h_s\le {\sqrt{8Q_0}\/\pi \r_n},
 \qqq M_n'\le {M_n\/\r_n},\qqq M_n''\le {M_n\/\r_n^2},
\]
\[
\lb{T5-2}
\sum_{n\in \Z} M_n\le {\pi ^2\/3\r^{2}}Q_0,\qq where \qqq \r=\inf_{n\in\Z}\r_n, 
\]
\[
\lb{T5-3}
{1 \/\r}\le {5\/2}e^{5\sqrt{Q_0/2}}.
\]
If $q'\in L^2(0,1)$ or $q_D'\in L^2(0,1)$, then
\[
\lb{T5-4}
 V_n(x)\le {2\/|x|\sqrt{|z_n^+z_n^-|}}
\rt(Q_0+{Q_2\/\r_n^2}\rt),\qqq  \ x\in \ol g_n.
\]
\end{theorem} 

\begin{corollary} \label{T6} 
Let real $q_D\in L^2(0,1)$. 
Then the following estimates hold

\[
\lb{T6-1}
\rt(\sum_{n\in\Z} 
|a_n-{|g_n|h_n\/4}|^{1\/3}\rt)^3 \le 
  Q_0^3
\rt({\pi ^2\/4^{7}\r^{4}}\rt)\rt( 1+{\sqrt{8Q_0}\/\pi \r}  \rt),
\]
\[
\lb{T6-2}
\sum_{n\in\Z} 
|a_n-{|g_n|h_n\/4}|\le 
{g_+^4\/2^7}{\pi ^2Q_0\/3\r^{4}}\rt(1+6{\sqrt{8Q_0}\/\pi \r}\rt),\qqq 
g_+=\sup_n |g_n|.
\]
\end{corollary}

In order to prove these results we consider the quasimomentum of the Hill (or the Dirac) operator  as a conformal mapping.
That makes possible to reformulate the problem for the differential operator as a problem of the conformal mapping theory (see \cite{KK1},
\cite{K1}-\cite{K8} and \cite{MO1}). We use  the Poisson integral  for the domain $\C_+\cup \C_-\cup g_n $ and some a priori estimates  from \cite{KK1}. The results of this paper are used in \cite{KK3}.

\section{Estimates for conformal mappings}
\setcounter{equation}{0}

Consider a conformal mapping $z:\cK_+\to\C_+$ with asymptotics
$z(iv)=iv(1+o(1))$ as $v\to\iy$, where  $z(k), k=u+iv\in \K$.
The domain $\cK_+=\C_+\cap \cK $ for some sequence $(h_n)_{n\in\Z}\in\ell^{\iy}, h_n\ge 0$ and the domain 
$\cK=\C \sm \cup_{n\in\Z}\G_n$, where $\G_n=[u_n -ih_n,u_n+ih_n]$ and 
$u_n, n\in \Z$ is strongly
increasing sequence of real numbers such that $u_n\to \pm \iy $ as
$n\to \pm \iy $. 
The difference of any two such mappings equals a real constant. Thus the imaginary part $y(k)=\Im z(k)$ is unique.
We call such mapping $z(k) $ the comb mapping.
Define the inverse mapping $k(\cdot): {\C _+}\to \cK_+$.
It is clear that $k(z), z=x+iy\in \C _+$ has the continuous extension
into $\ol {\C _+} $. We define "gaps" $g_n$, "bands"  $s_n$ and the "spectrum" $s$ of the comb mapping by:
$$
g_n=(z_n^-, z_n^+)=(z(u_n-0), z(u_n+0)),\qqq
s_n=[z_{n-1}^+, z_n^-],\qqq s=\cup_{n\in\Z }s_n.
$$
It is well known that the set $s$ can not be the spectrum of two different comb mappings \cite{Le}.
The function $u(z)=\Re k(z) $ is strongly  increasing on each band
$s_n$ and $u(z)=u_n$  for all $z\in [z_n^-,z_n^+],\ n\in\Z $;
the function $v(z)=\Im k(z)$ equals zero  on each  band $s_n$ and
is  strongly convex on each gap $g_n\ne \es$ and has the maximum  at some point $z_n $ given by $v(z_n)=h_n$.
 If  the gap  is empty we set $z_n=z_n^{\pm}$. The function $z(\cdot) $ has an analytic extension (by the symmetry) from  the domain $\cK_+$   onto  the domain $\cK$
and $z(\cdot): \cK\to z(\cK)=\cZ=\C\sm\cup \ol g_n $ is  a conformal
mapping. 
For each gap $g_n=(z_n^-,z_n^+)\ne \es$  we define the effective masses $m_n^{\pm}$ by
\[
\lb{2.20}
z(k)-z_n^{\pm} ={(k-u_n)^2 \/  2m_n^{\pm}}+ O((k-u_n)^3)
\qqq {\rm as} \ \ z\to z_n^{\pm}.
\]
If $|g_n|=0$, then we set $m_n^{\pm}=0$. 
It is clear that if  $|g_n|>0$, then $\pm m_n^{\pm}>0$.
 These and  others properties of the comb mappings it is possible  to  find in the papers of  Levin \cite{Le}.

We emphasize that the comb mapping $k(z)$ with
$u_n=\pi n, n\in \Z$  corresponds to the quasimomentum for the
Dirac operator with periodic coefficients.
If the comb mapping $k(z)$ is symmetric, $k(-z)=-k(z), z\in \cZ$,
with $u_n=\pi n, n\in \Z$, then it
corresponds to the  quasimomentum  for  the   Hill operator. Furthermore, the comb mapping  is  the rotation number for
the Schr\"odinger operator with  some  almost periodic  potential
(see [JM]).

For some $(h_n)_{n\in \Z}\in \ell_{\R }^{\iy}$ we fix some gap $G=(-c,c), c>0$.
The function $v=\Im k(z)$ satisfies
\[
\lb{6.9}
v(x)=r(x)\big(1+V(x)\big),\qq
V(x)={1\/\pi}\int_{\R \sm g}{v(t)dt\/|t-x|\,r(t)}, \qq 
r(x)=|c^2-x^2|^{1\/2},\ 
x\in (-c,c),
\]
see [KK1]. Recall that  $v''(x)<0$ for all $x\in G$. The function $v(x), x\in G$ has a unique maximum $h$ at some point  $z_0\in (-c,c)$ given by
\[
\lb{ihg}
 h=v(z_0)=r(z_0)(1+V_0)\ge c,\qqq V_0=V(z_0).
\]

\begin{lemma}
\lb{T6.4}
Let a comformal mapping $k:\cZ \to \cK$
have a gap $G=(-c,c),c>0$ for some $(h_n)_{n\in \Z}\in \ell^\iy$ and 
$V_0=V(z_0),\ V_0'=V'(z_0),.., r_0=r(z_0),\n=-1/v''(z_0)>0$. Then
\[
\lb{6.18}
z_0={r_0^2V_0'\/1+V_0},\qqq where \qq  r_0\le c,\ 
\]
\[
\lb{6.19}
|c-r_0|\le {r_0^3{V_0'}^2\/2(1+V_0)^2},
\]
\[
\lb{6.20}
|h-c(1+V_0)|\le {r_0^3V_0'^2\/2(1+V_0)},
\]
\[
\lb{6.21}
|{c^2\/h}-\n|\le h c^2((3+c){V_0'}^2+V_0'').
\]
\end{lemma}
\no {\bf Proof.} We have $v'(z_0)=0$ and
$0=v'(x)=r'(x)(1+V(x))+r(x)V'(x)$ at $x=z_0$ and 
$(r^2)'=-2x$, which yield \er{6.18}.
Substituting \er{6.18} into the identity $z_0^2=c^2-r_0^2$, where $r_0=r(z_0)$ we obtain
$
c^2-r_0^2={r_0^4V_0'^2\/(1+V_0)^2},
$
which implies \er{6.19}.

Estimate \er{6.19} together with the identity 
$h=r_0(1+V_0)$ yield
$$
|h-c(1+V_0)|=(1+V_0)|r_0-c|\le {r_0^3V_0'^2 \/2(1+V_0)},
$$
which gives \er{6.20}.
Differentiating the function $v'(x)r(x)$ at $z_0$ we obtain
$$
r_0v''(z_0)= (rv')'|_{z=z_0}=\lt(-x(1+V)+r^2
V'\rt)'|_{z=z_0}=-(1+V_0)-3z_0V_0'+r_0^2V_n'',
$$
which gives
$
r_0v''(z_0)=-(1+V_0)-3z_0V_0'+r_0^2V_0''.
$
Then the definition $\n_n=-1/v''(z_0)$ implies
$r_0=\n (1+V_0+3V_0' z_0-r_0^2V_0'')$, which yields
\[
{r_0\/1+V_0}={r_0^2\/h}=\n+R,\ \ R=\n {3V_0' z_0-r_0^2V_0''\/1+V_0},\qqq
c^2-r_0^2={r_0^4V_0'^2\/(1+V_0)^2}.
\]
Thus ${c^2\/h}=\n+R_1, R_1=R+{c^2-r_0^2\/h}$. Thus
$|R|\le \n r_0^2(3 {V_0'}^2+V_0'')$ and 
${c^2-r_0^2\/h}\le {r_0^3V_0'^2\/(1+V_0)^2}$, which 
yield \er{6.21}, where the estimate $\n\le h$ from \cite{KK4} was used.
$ \BBox $

Recall the identity from \cite{KK1}
\[
\lb{iemg}
\pm m^{\pm}=c(1+V(\pm c)^2\ge c.
\]
We obtain estimates of the action.
Let $M=\max_{x\in G} V(x), M'=\max_{x\in G} V'(x),...$

\begin{lemma}
\label{T6.5}
Let a comformal mapping $k:\cZ \to \cK$
have a gap $G=(-c,c),c>0$ for some $(h_n)_{n\in \Z}\in \ell^\iy$. Then the actions $a={1\/\pi}\int_{G}v(x)dx$
and $A={4\/\pi}\int_{G}xv(x)dx$ and $b=|a|^{1\/2}>0$ satisfy
\[
\lb{6.25}
|a-{c^2\/2}(1+V(0))|\le {c^4\/8} M'',\qqq 
\]
\[
\lb{leA}
|A|\le {c^4\/2}M', \qqq 
\]
\[
\lb{6.26}
|a-{ch\/2}|\le {c^4\/8}(M''+6{M'}^2),
\]
\[
\lb{am1}
|a-{c\/2}\sqrt{c|m_\pm|}|\le {c^3\/2} (M'+{c\/4} M''),
\]
\[
\lb{6.27}
|b-c\sqrt{(1+V(0))/2}|\le { c^{3}\/8\sqrt2} M'',
\]
\[
\lb{6.28}
|b-\sqrt{ch/2}|\le { c^3\/8\sqrt 2}(M''+6{M'}^2),
\]
\[
\lb{am2}
|b-{c^{3/4}\/\sqrt2}|m_\pm|^{1\/4}|\le {c^2\/2\sqrt2}( M'+{c\/4} M'').
\]
\end{lemma}
\no {\bf Proof.} Using the Taylor formula 
$
V(x)= V(0)+ V'(0)x+\int _0^{x} V''(s)(x-s)ds,
$
we obtain
$$
\pi a=\int_{G}v(x)dx=\int_{G}r(x)(1+V(x))dx=a_{(1)}+a_{(2)},\qq
a_{(2)}=\int_{G}\!\!r(x)dx\int_{0}^{x}\!\! V''(s)(x-s)ds,
$$
$$
a_{(1)}=\int _{G}r(x)(1+V(0)+V'(0)x)dx=\pi {c^2\/2}(1+V(0)).
$$
Consider $a_{(2)}$. The function $V''(x)>0$ for all $x\in G$,
then for some $x_1\in G$ we get
$$
a_{(2)}=\int_{G}r(x)dx \int _{0}^{x} V''(s)(x-s)ds=
V''(x_1)\int _{G}r(x)dx \int _{0}^{x}(x-s)ds=
\pi V''(x_1){c^4\/2^4}
$$
which yields \er{6.25}.
Using $V(x)= V(0)+ \int _0^{x} V'(s)ds$ and $\int_{G}xr(x)(1+V(0))dx=0$, we have
$$
A={4\/\pi}\int_{G}xv(x)dx={4\/\pi}\int_{G}xr(x)dx\int _0^{x} V'(s)ds,
$$$$
|A|\le{4\/\pi}\int_{G}|x|r(x)|\int _{0}^{x} V'(s)ds|dx\le 
{4M'\/\pi}\int_{G}x^2r(x)dx={c^4M'\/ 2},
$$
which yields \er{leA}. Using \er{ihg} we have the following identity 
$$
h=c(1+V(0))+(h-c(1+V_0)) +cz_0V'(x_1), \qq where \qq V_0=V(z_0),
$$
for some $x_1\in G$, which yields 
$$
a-{hc\/2}=\rt(a-{c^2\/2}(1+V(0))\rt)
+{c(h-c(1+V_0))\/2}+{c^2z_0\/2}V'(x_1).
$$
Using \er{6.25}, \er{6.20}, \er{6.18}  we obtain
$$
|a-{hc\/2}|\le {c^4\/2^3} M''+{c^4M'^2\/4}+{c^4{M'}^2\/2}<
{c^4\/8}(M''+6{M'}^2),
$$
which yields \er{6.26}. Using \er{iemg} we obtain
$$
|a-{c\/2}\sqrt{cm_\pm}|\le |a-(c^2/2)(1+V(0))|+{c^{2}\/2}
|V(0)-V(\pm c)|\le {c^4\/8} M''+{c^3\/2} M',
$$
which yields \er{am1}. Using \er{6.25}  we obtain
$
|b-c\sqrt{(1+V(0))/2}|\le { c^4\/8} {M''\/|b+2^{-{1\/2}}c|},
$
and the estimate $a=b^2\ge c^2/2$
 yields \er{6.27}. Using \er{6.26}  we obtain
$
|b-\sqrt{ch/2}|\le { c^4\/8} {M''+6{M'}^2\/|b+\sqrt{ch/2}|},
$
and the estimate $a=b^2\ge c^2/2, h\ge c$
 yields \er{6.28}. Due to \er{am1}  we have
$$
|b-{c^{3/4}\/\sqrt2}|m_\pm|^{1\/4}|\le {c^3\/2}{( M'+{c\/4} M'')\/(b+
{c^{3/4}\/\sqrt2}|m_\pm|^{1\/4})},
$$
and $|m_\pm|\ge c, a\ge {c^2\/2}$ gives \er{am2}.
\BBox

\begin{lemma} \label{Tefm}
Let a comformal mapping $k:\cZ \to \cK$
have a gap $G=(-c,c),c>0$ for some $(h_n)_{n\in \Z}\in \ell^\iy$.  Then 
\[
\lb{h+}
|h-h^{\pm }|\le 2c^2\rt(M'+
{{V_0'}^2\sqrt{cm^\pm}\/4(1+V_0)^2}\rt), \qqq
where \ h^{\pm}=\sqrt{c|m^{\pm}|},
\]
\[
\lb{m+-}
|m^++m^-|\le 4c^2(1+M)M',
\]
\[
\lb{m+-0}
|m^{\pm}\mp m|\le 4c^2(1+M)\rt(M'+
{{V_0'}^2\sqrt{cm^\pm}\/4(1+V_0)^2}\rt),\qqq
where \ \ m={h^2\/c}.
\]
\end{lemma}
\no {\bf Proof.}  
Let $V_0=V(z_0)$ and $V_\pm=V(\pm c)$.  Identities \er{iemg}, \er{ihg} give
$$
|h-h^\pm|=|r_0(1+V_0)-c(1+V_\pm)|\le r_0|V_0-V_\pm| +|r_0-c|(1+V_\pm).
$$
Then estimates \er{6.19}, \er{iemg} give 
$|r_0-c|(1+V_\pm)\le 
{c^2{V_0'}^2\sqrt{c|m^\pm|}\/2(1+V_0)^2},
$
and we also have $|V(z_0)-V(\pm c)|\le 2cM'$. Combining two last estimates we obtain \er{h+}.

We prove \er{m+-}. Identity \er{iemg} gives
$$
|m^++m^-|=c(2+V_++V_-)|V_+-V_-|\le
4c^2(1+M)M',
$$
which yields \er{m+-}. 
We prove \er{m+-0}. Using \er{iemg}, \er{ihg} we have
$$
|m^+-m|=
{|(h^+)^2-h^2|\/c}=
{|h^+-h||h^++h|\/c}
\le
|h^++h|2c\rt(M'+
{{V_0'}^2\sqrt{cm^\pm}\/4(1+V_0)^2}\rt),
$$
and $h+h^+\le 2c(1+M)$ yields \er{m+-0}. \BBox

We rewrite \er{6.9} for any gap $g_n=(z_n^-,z_n^+)$ by
\[
v(x)=v_n(x)(1+V_n(x)),\qq
V_n(x)={1\/\pi}\int_{g\sm g_n}{v(t)dt\/v_n(t)|t-x|}, \ 
v_n(x)=|(x-z_n^-)(x-z_n^+)|^{1\/2},
\]
$ x\in g_n$.
We have the following identities and estimates
\[
\lb{6.11}
V_n'(x)={1\/\pi} \int_{g \sm g_n}{v(t)\,dt\/(t-x)|t-x|v_n(t)},\qq
V_n''(x)={2\/\pi}\int_{g\sm g_n}{v(t)\,dt\/|t-x|^3\,v_n(t)}>0,
\]
\[
\lb{6.13}
|V_n'(x)|\le {V_n(x)\/\r_n},\qq \qq
0<V_n''(x)\le {2V_n(x)\/\r_n^2},
\]
for all $x\in g_n$.
Define the numbers $M_n=\max_{x\in g_n} V_n(x), M_n'=\max_{x\in g_n} V_n'(x),$ . Let $Q_p={1\/\pi}\int_\R x^{2p}v(x)dx, p\ge 0$.
We prove Theorem \ref{T5} for more general case.

\begin{lemma} \label{TeM} Let $u_*=\inf_{n\in \Z} (u_{n+1}-u_n)>0$
and let $Q_0<\iy$.  Then 
\[
\lb{TeM-1}
 M_n\le {2\/\pi \r_n}\sup_{s\in \Z} h_s\le {\sqrt{8Q_0}\/\pi \r_n}, \qqq M_n'\le {M_n\/\r_n},\qqq M_n''\le {M_n\/\r_n^2},
\]
\[
\lb{TeM-2}
\sum M_n\le {\pi ^2\/3\r^{2}}Q_0,
\]
\[
\lb{TeM-3}
\r=\inf \r_n\ge {2u_*\/5\pi}\exp (-{5\pi\/2u_*}\sup_{s\in \Z} h_s).
\]
If additionly $Q_2<\iy$, then 
\[
\lb{TeM-4}
 V_n(x)\le {2\/|x|\sqrt{|z_n^+z_n^-|}}
\rt(Q_0+{Q_2\/\r_n^2}\rt),\qqq  \ x\in \ol g_n.
\]
\end{lemma}
\no {\bf Proof.} Let $h_+=\sup h_n$. We have for $x\in g_n$ 
$$
V_n(x)={1\/\pi}\int_{g\sm g_n}{v(t)dt\/|t-x|v_n(t)}
\le {h_+\/\pi}\int_{g\sm g_n}{dt\/|t-x|v_n(t)}\le
{h_+\/\pi}\int_{|s|>\r_n}{ds\/|s|^2}\le{2h_+\/\pi\r_n}
$$
and using $h_+^2\le 2Q_0$ (see \cite{K5}) we obtain the first estimate in \er{TeM-1}.  The two last estimates follow from \er{6.13}.

Let $M_n=V_n(x_n)$ for some  $x_n\in g_n$. Using $|t-x_n|\ge \r |m-n|$
for any $t\in g_m$, we obtain
$$
\sum M_n\le \sum_n\sum _{m\ne n}{1\/\pi}\int_{g_m}
{v(t)dt\/|t-x_n|v_n(t)}\le 
\sum_n\sum _{m\ne n}{1\/\pi}\int_{g_m}
{v(t)dt\/\r^2|n-m|^2}={Q_0C\/\r^2},
$$
where $C=2\sum_{n\ge 1}n^{-2}=2{\pi^2\/6} ={\pi^2\/3}$.
Estimate \er{TeM-3} was proved in \cite{KK4}.

 Using the simple estimate ${1\/|t-z|}\le {1\/|z|}\rt(1+
{|t|\/\r_n}  \rt)$ for all $z\in \ol g_n, t\in g\sm g_n$, we obtain
$$
V_n(x)={1\/\pi}\int_{g \sm g_n}{v(t)dt\/|t-x|v_n(t)}
\le {1\/\pi}\int_{g \sm g_n}\rt(1+{|t|\/\r_n}\rt)^2{v(t)dt\/|x|\sqrt{|z_n^+z_n^-|}}
$$
$$
\le  {2\/\pi |x|\sqrt{|z^+z^-|}}\int_{\R}\rt(1+{t^2\/\r_n^2}\rt)v(t)dt\le {2\/|x|\sqrt{|z_n^+z_n^-|}}\rt(Q_0+{Q_2\/\r_n^2}\rt),
$$
which yields \er{TeM-4}.
\BBox

\section { Estimates for the Hill and the Dirac   operators}
\setcounter{equation}{0}

Firstly, we will apply the results from Section 2 for
the case the Dirac operator $T_D$.

 {\bf Proof of Theorem \ref{TD1}.}
 Estimate \er{TD1-1}, \er{TD1-2}, \er{TD1-3},\er{TD1-4},\er{TD1-5} follows from \er{6.18}, \er{6.20}, \er{h+}, \er{m+-}, \er{m+-0}
 respectively.
 \BBox

{\bf Proof of Theorem \ref{TD2}.}
Estimate \er{TD2-1}, \er{TD2-2}, \er{TD2-3}, \er{TD2-4} follows from \er{6.25}, \er{6.26}, \er{am1}, \er{6.27}.
\BBox


Secondly we will apply the results from Section 2 for
the case the Hill operator $T= -{d^2\/dt^2}+(q_0+q'(t))$ in $L^2(\R)$ with  the 1-periodic real potential $q\in L^2(0,1)$. 
Recall that $\l_0^+=0$ for the Hill operator. In this case we have estimates
\[
{|\g_n|\/ 2\sqrt {\l_n^+}} \le |g_n|={|\g_n|\/z_n^++z_n^-}
\le {|\g_n|\/2n\r},\ \ m_n^{\pm}=2z_n^{\pm}\m_n^{\pm},\ \ \ n\ge 1.
\]

 {\bf Proof of Theorem \ref{T1}.}
Recall $\l_n^0={\l_n^++\l_n^-\/2}$ and $\l_n=z_n^2$. Using the identity
 $$
4(z_n^0)^2=(z_n^-+z_n^+)^2=\l_n^-+\l_n^++2z_n^-z_n^+=
4\l_n^0-|g_n|^2,
$$
and the estimate \er{6.18} we obtain
$$
|\l_n^0-{|g_n|^2\/4}-\l_n|=|(z_n^0)^2-\l_n|=
|z_n^0+z_n||z_n^0-z_n|\le
|z_n^0+z_n|{|g_n|^2\/4}|V_n'(z_n)|
\le {3|\g_n|^2\/8z_n^0}|V_n'(z_n)|,
$$
which yields \er{T1-1}. Recall $2z_n^\pm\m_n^\pm=m_n^\pm$.
Using \er{6.20} we obtain \er{T1-2}.

Using \er{h+} and $|g_n||m_n^\pm|\le 2|\g_n||\m_n^\pm|$ we obtain \er{T1-3}.

Using \er{m+-}, \er{iemg} and $2z_n^\pm\m_n^\pm=m_n^\pm$ we obtain
$$
|\m_n^++\m_n^-|={|z_n^0(m_n^++m_n^-)+c_n(m_n^--m_n^+)|\/2z_n^+z_n^-}\le
{z_n^0|g_n|^2\/2z_n^+z_n^-}(1+M_n)M_n'+{|g_n|^2\/4z_n^+z_n^-}(1+M_n),
$$
which yields \er{T1-4}. 
 Using $\m_n={m_n\/2z_n^0}={2h_n^2\/|\g_n|}$ we obtain 
$
\m_n^\pm\mp \m_n={\mp |g_n| m_n^\pm\/4z_n^0z_n^\pm}+{m_n^\pm-m_n\/2z_n^0}
$.
The estimate \er{m+-0} yields (note that $|g_n||m_n^\pm|\le 2|\g_n||\m_n^\pm|$)
$$
{|m_n^\pm-m_n|}\le {|g_n|^2}
(1+M_n)M_n'\rt(1+{M'_n\/8}\sqrt{2|g_n||m_n^\pm|}\rt)
$$
and \er{iemg} gives $| m_n^\pm|\le |g_n|(1+M_n)^2/2$. The two last
estimates imply \er{T1-5}.
 \BBox

{\bf Proof of Theorem \ref{T2}.}  Using \er{leA} we obtain
\[
\lb{32}
A_n={4\/\pi}\int_{g_n}xv(x)dx=4z_n^0a_n+A_{n0},\qq
A_{n0}={4\/\pi}\int_{g_n}(x-z_n^0)v(x)dx ,\qq
|A_{n0}|\le {|g_n|^4\/2^5}M_n'.
\]
Estimate \er{6.25} and $2z_n^0|g_n|=|\g_n|$  imply
\[
\lb{A1}
|4z_n^0a_n-(|\g_n|^2/8z_n^0)(1+V_n(z_n^0))|\le z_n^0{|g_n|^4\/2^5} M_n'',
\]
which together with \er{32} yields \er{T2-1}.
Estimate \er{6.26} implies
$$
|4z_n^0a_n-{|\g_n|h_n\/2}|\le z_n^0{|g_n|^4\/2^5}(M_n''+6{M_n'}^2),
$$
which together with \er{32} yields \er{T2-2}.
Estimate \er{am1} implies
$$
|4z_n^0a_n-{|\g_n|\/2}|c_nm_n^\pm|^{1\/2}|\le 2z_n^0
c_n^3 (M_n'+c_nM_n'')
$$
which together with \er{32} yields \er{T2-3}.
Using and $\int_{g_n}(x-z_n^0)v_n(x)dx=0$, we obtain
\[
\lb{eb}
\b_n^2=A_n={4\/\pi}\int_{g_n}xv(x)dx\ge {4\/\pi}\int_{g_n}xv_n(x)dx=
{4z_n^0\/\pi}\int_{g_n}v_n(x)dx={z_n^0|g_n|^2\/2}
\]
Recall that \er{T2-1} yields
$
\rt|A_n-|\g_n|^2{(1+V_n(z_n^0))\/8z_n^0}\rt|\le {|g_n|^4\/2^5}(z_n^0 M_n''+M_n'),
$
which gives
$$
\rt|\b_n-|\g_n|\sqrt{(1+V_n(z_n^0))\/8z_n^0}\rt|\le {|g_n|^4\/2^5
(\b_n+|g_n|\sqrt{z_n^0/2})}(z_n^0 M_n''+M_n'),
$$
and \er{eb}  yields \er{T2-4}. 
Using \er{T2-2},\er{T2-3} and the similar arguments we obtain
\er{T2-5},\er{T2-6}. 
 \BBox

{\bf Proof of Corollary \ref{T6}}.  
Using \er{T5-1}, \er{T5-2} we obtain

$$
\sum_{p\in\Z} M_p''\le {1\/\r^2}\sum_{p\in\Z} M_p\le 
{\pi ^2Q_0\/3\r^{4}},\qqq \sum_{p\in\Z} {M_p'}^2\le 
 {\sqrt{8Q_0}\/\pi \r} {\pi ^2Q_0\/3\r^{4}},
$$
which yields
\[
\lb{c1}
\sum_{p\in\Z} (M_p''+6{M_p'}^2)\le 
{\pi ^2Q_0\/3\r^{4}}(1+6{\sqrt{8Q_0}\/\pi \r}).
\]
Using \er{TD2-2} and \er{c1}  we obtain
$$
\sum_{n\in\Z} 
|a_n-{|g_n|h_n\/4}|^{1\/3}\le \sum_{n\in\Z} {|g_n|^{4\/3}\/2^{7\/3}}(M_n''+6{M_n'}^2)^{1\/3}\le 2^{-7\/3}
\rt(\sum_{n\in\Z}|g_n|^2\rt)^{2\/3}  \rt(\sum_{p\in\Z} (M_p''+6{M_p'}^2)  \rt)^{1\/3}
$$
$$
\le 2^{-7\/3}\rt({Q_0\/8}\rt)^{2\/3}
\rt({\pi ^2Q_0\/3\r^{4}}\rt)^{1\/3}\rt( 1+{\sqrt{8Q_0}\/\pi \r}  \rt)^{1\/3},
$$
which yields \er{T6-1}.
Using \er{TD2-2} and \er{c1}  we obtain
$$
\sum_{p\in\Z}|a_n-{|g_n|h_n\/4}|\le \sum_{p\in\Z}{|g_n|^4\/2^7}(M_n''+6{M_n'}^2)\le
{g_+^4\/2^7}{\pi ^2Q_0\/3\r^{4}}(1+6{\sqrt{8Q_0}\/\pi \r}),
\qq g_+=\sup |g_n|.
$$

\BBox

\end{document}